\documentclass{article}

\def\ss{\subseteq}

\def\dbar{\overline{\partial}}

 \def\HollowBox #1#2{{\dimen0=#1 \advance\dimen0 by -#2       
       \dimen1=#1 \advance\dimen1 by #2                       
        \vrule height #1 depth #2 width #2                    
        \vrule height 0pt depth #2 width #1                   
        \llap{\vrule height #1 depth -\dimen0 width \dimen1}%
       \hskip -#2                                             
       \vrule height #1 depth #2 width #2}}                   
 \def\BoxOpTwo{\mathord{\HollowBox{6pt}{.4pt}}\;}             

\def\endpf{\hfill $\BoxOpTwo$}

\font\teneufm=eufm10
\font\seveneufm=eufm7
\font\fiveeufm=eufm5
\newfam\eufmfam
\textfont\eufmfam=\teneufm
\scriptfont\eufmfam=\seveneufm
\scriptscriptfont\eufmfam=\fiveeufm

\newfam\msbfam
\font\tenmsb=msbm10  \textfont\msbfam=\tenmsb
\font\sevenmsb=msbm7  \scriptfont\msbfam=\sevenmsb
\font\fivemsb=msbm5    \scriptscriptfont\msbfam=\fivemsb
\def\Bbb{\fam\msbfam \tenmsb}

\def\CC{{\Bbb C}}

\newtheorem{theorem}{Theorem}

\newtheorem{remark}[theorem]{Remark}

\makeindex

\usepackage{graphicx}

\usepackage{amsmath}

\begin{document}

\begin{center}
\huge \bf
The Corona Problem with Two Pieces of Data
\end{center}
\vspace*{.12in}

\begin{center}
\large Steven G. Krantz\footnote{Author supported in part
by the National Science Foundation and by the Dean of the Graduate
School at Washington University.}\footnote{{\bf Key Words:}  corona problem,
bounded holomorphic functions, domain of holomorphy, ball}\footnote{{\bf MR Classification
Numbers:}  30H80, 32A38, 32A65}
\end{center}
\vspace*{.15in}

\begin{center}
\today
\end{center}
\vspace*{.2in}

\begin{quotation}
{\bf Abstract:} \sl
We study the corona problem on the unit ball in $\CC^n$, and
more generally on strongly pseudoconvex domains in $\CC^n$.  When
the corona problem has just two pieces of data, and an extra geometric
hypothesis is satisfied, then we are able to
solve it.
\end{quotation}
\vspace*{.25in}

\setcounter{section}{-1}

\section{Introduction}

Banach algebras were created in the thesis of I. M. Gelfand in 1938.  Some
of the key ideas in the subject are recorded in [GEL1] and [GEL2].  In 1942,
S. Kakutani formulated the following fundamental question (known as
the {\it corona problem}) in the subject:
\begin{quote}
Suppose that $f_1, \dots, f_k$ are bounded, holomorphic
functions on the unit disc $D$ with the property that
$$
|f_1(\zeta)| + |f_2(\zeta)| + \cdots + |f_k(\zeta)| > \delta > 0
$$
for some positive constant $\delta$ and all $\zeta \in D$.  
[We call $f_1, \dots, f_k$ a set of {\it corona data} on the disc.] 
Do there exist bounded, holomorphic $g_j$, $j = 1, \dots, k$, on $D$ such that
$$
f_1(\zeta) \cdot g_1(\zeta) + f_2(\zeta) \cdot g_2(\zeta) + \cdots +
       f_k(\zeta) \cdot g_k(\zeta) \equiv 1
$$
for all $\zeta \in D$?
\end{quote}

It turns out that this question is equivalent to asking whether the
point evaluation functionals are weak-$*$ dense in the space
of multiplicative linear functionals on the Banach algebra
$H^\infty(D)$.  This equivalence is explained in some detail
in [KRA2].  An authoritative treatment of the corona problem on the
disc appears in [GAR].  

Lennart Carleson [CAR] solved the corona problem on the disc $D$ in the affirmative
in 1962.  Since that time, the corona problem has been proved to
have a positive solution on a variety of different types of domains in $\CC$.
See, for example [GAJ].  Brian Cole produced a Riemann surface on which the corona
problem fails (see [GAM] and [GAR] for discussion).  

Meanwhile, there have been investigations of the corona problem on domains in $\CC^n$, particularly
on the ball and the polydisc.  A variety of counterexamples to the corona problem in
several variables have been produced in [SIB] and [FOS].  

At this time there is no known domain in the plane $\CC$ on
which the corona problem is known to fail. And there is no
known domain in $\CC^n$, for $n \geq 2$, on which the corona
problem is known to hold true. Some partial results on the
corona problem in several complex variables appear, for
instance, in [KRL1], [KRL2], [KRL3], and [BCL].

\section{Fundamental Results}

In the present paper we study the corona problem on
the ball in $\CC^n$, and more generally on strongly
pseudoconvex domains in $\CC^n$, when there are
just two pieces of corona data (i.e., $k = 2$).  Using
the Koszul complex (an idea developed for the computation
of the cohomology of a Lie algebra), we are able
to produce a solution of the corona problem in these circumstances.
Afterwards we shall make some remarks about the case when
there are an arbitrary (but, of course, finite) number of
pieces of corona data.

We shall work on a {\it domain} $\Omega$ in $\CC^n$, where here
a domain is a connected, open set.

We take this opportunity to review how the Koszul complex works
in the present context.  Refer, for instance, to [KRA1] or [KRA2] for
more information about the Koszul complex.

If $f_1, f_2, \dots, f_k$ are given corona data on a domain $\Omega \ss \CC^n$,
we define 
$$
U_j = \left \{\zeta \in \Omega: |f_j(\zeta)| > \frac{\delta}{2k} \right \} \, .
$$
It is easy to see then that the $U_j$ are open and $\cup_j U_j \supseteq \Omega$.  Let $\{\varphi_j\}$ be 
a partition of unity subordinate to the covering $\{U_j\}$.  
					   
Now set
$$
g_j = \frac{\varphi_j}{f_j} + \sum v_{ji} f_i \, ,  \eqno (\star)
$$
where the $v_{ji}$ are functions yet to be determined and $v_{ij} = - v_{ji}$ and $v_{jj} \equiv 0$.  Then it is a simple
formal verification to see that the $g_j$ are well-defined and
$$
\sum_j f_j g_j \equiv 1 \, .
$$

Now our main result is this:

\begin{theorem} \sl
Let $\Omega \ss \CC^n$ be strongly pseudoconvex.  Let $f_1, f_2$ be
corona data on $\Omega$.  Let ${\cal Z}_1 = \{z \in \Omega: f_1(z) = 0\}$
and ${\cal Z}_2 = \{z \in \Omega: f_2(z) = 0\}$.  Assume that
there is an $\eta > 0$ such that, if $z \in {\cal Z}_1$ and $w \in {\cal Z}_2$ then
$|z - w| > \eta$.  Then there exist bounded, holomorphic $g_1, \dots, g_k$ such
that 
$$
f_1 g_1 + f_2 g_2 \equiv 1 \, .
$$
\end{theorem}
\noindent {\bf Proof:}  Let us return to the Koszul complex.
Because the $\varphi_j$ are real functions, we see immediately that
the $g_j$ so defined are not necessarily holomorphic.  But we would
like to see that it is possible to choose $v_{ji}$ so that each
$g_j$ will in fact be holomorphic.  We may assume that $\varphi_1 \equiv 1$ on $U_1 \setminus U_2$ and
$\varphi_2 \equiv 1$ on $U_2 \setminus U_1$.This entails
$$
\dbar g_j = \frac{\dbar \varphi_j}{f_j} + \sum (\dbar v_{ji}) f_i \equiv 0 \, .
$$
In general the algebra involved in understanding this last equation
is rather complicated.  But we take $k = 2$.  Then we find the situation
to be tenable.  Namely, we may write
$$
\dbar v_{11} f_1 + \dbar v_{12} f_2 = - \frac{\dbar \varphi_1}{f_1}
$$
and
$$
\dbar v_{21} f_1 + \dbar v_{22} f_2 = - \frac{\dbar \varphi_2}{f_2} \, .
$$
This simplifies to
$$
\dbar v_{12} f_2 = - \frac{\dbar \varphi_1}{f_1}
$$
and
$$
\dbar v_{21} f_1 = - \frac{\dbar \varphi_2}{f_2} \, .
$$

Using the given identities on the $v_{ji}$, we then find that
$$
\dbar v_{12} f_2 = - \frac{\dbar \varphi_1}{f_1}
$$
and
$$
- \dbar v_{12} f_1 = - \frac{\dbar \varphi_2}{f_2} \, .
$$
Since $\varphi_1 + \varphi_2 \equiv 1$, we know that $\dbar \varphi_1 = - \dbar \varphi_2$.
So in fact the last two equations say the same thing:
$$
\dbar v_{12} = \frac{\dbar \varphi_2}{f_1 f_2}  \equiv \lambda \, .   \eqno (*)
$$

This last is the $\dbar$-problem that we must solve.  First
notice that the righthand side is $\dbar$-closed.  So the
problem can be solved on a pseudoconvex domain (see [KRA1]).  

Now we have assumed that the zero set of $f_1$ is bounded
away from the zero set of $f_2$.  That is, there is an $\eta > 0$ so
that, if $f_1(z) = 0$
and $f_2(w) = 0$ then $|z - w| > \eta$.  [Unfortunately, the corona condition
does not guarantee this assumption.  It {\it does} guarantee that the
assumption holds on the smaller disc $\overline{D}(0,r)$, $0 < r < 1$, but
of course the estimate on $\lambda$ does not hold uniformly in $r$.]

As a result of this extra hypothesis, it is easy to arrange for $\dbar \varphi_1$
(resp.\ $\dbar \varphi_2$) to be nonzero only when $|f_1| \geq \delta/2k$ 
and $|f_2| \geq \delta/2k$.  

As a consequence, the data function on the righthand side of $(*)$ is bounded,
and we may solve to obtain a bounded $v_{12}$.  Therefore $g_1$ and $g_2$, according
to $(\star)$, are bounded functions.  And we have solved the corona problem
for the data $f_1$, $f_2$.
\endpf
\smallskip \\

\begin{remark} \rm
An essential part of the argument just presented is that the
zero sets of $f_1$ and $f_2$ are separated.  One may
wonder to what degree this sort of property can be guaranteed.  In particular,
is it possible for a bounded holomorphic $h$ on the ball $B$ to have
zero set accumulating at every boundary point (if so, two such functions
could {\it not} have their zero sets separated).

In dimension 1, the answer to this last question is easily ``yes''.  For
let $\Omega = D$ be the unit disc.  Set ${\cal S}$ be the set 
consisting of
\begin{itemize}
\item Two evenly spaced points at distance $2^{-2}$ from the boundary of $D$;
\item Three evenly spaced points at distance $2^{-3}$ from the boundary of $D$;
\item Four evenly spaced points at distance $2^{-4}$ from the boundary of $D$;
\item \ \ \ \dots
\item $k$ evenly spaced points at distance $2^{-k}$ from the boundary of $D$.
\end{itemize}
Then the Blaschke product with zero set ${\cal S}$ will have the desired property.

In higher dimensions it is more difficult to construct such a function, but it
can be done.  The paper [BER] gives the tools for producing an example on the ball $B$ in complex
dimension 2.  Indeed, Berndtsson shows that any analytic set of bounded volume
is the zero set of a bounded function.  Thus one may take $\{p_j\}$ to be a countable,
dense set in $\partial B$ and then let $\{{\bf d}_{jk}\}$ be analytic
discs of radius $2^{-j+k}$, $k = 1, 2, \dots$, which lie in the domain and approach $p_j$.  Then the union
of all these discs is an analytic set with finite 2-dimensional volume, and is thus
the zero set of an $H^\infty$ function; that function has the properties that we seek.
Also note that the paper [HAS] gives an example on $B$ in all dimensions; indeed Hakim and Sibony
produce a bounded, holomorphic function whose zeros accumulate (but {\it not} admissibly!)
at each boundary point, but which is not the identically zero function.
\end{remark}

\section{Concluding Remarks}

The corona problems is an important part of classical function theory, and of
modern harmonic analysis.  The present paper offers the first positive
result in several complex variables that is in the spirit of the original
corona problem.  We hope in future papers to extend the result
to an arbitrary number of pieces of corona data.

\newpage

\noindent {\Large \sc References}
\bigskip  \\

\begin{enumerate}						    
	
\item[{\bf [BER]}] B. Berndtsson, Integral formulas for the
$\dbar$-equation and zeros of bounded holomorphic functions in
the unit ball, {\it Math.\ Ann.} 249(1980), 163--176.

\item[{\bf [BCL]}] B. Berndtsson, S. Y. Chang, and K. C. Lin,
Interpolating sequences in the polydisc, {\it Trans.\ Amer.\
Math.\ Soc.} 302(1987), 161--169.
								   
\item[{\bf [CAR]}] L. Carleson, Interpolation by bounded
analytic functions and the corona problem, {\it Ann.\ of
Math.} 76(1962), 547--559.
										      
\item[{\bf [FOS]}] J. E. Forn\ae ss and N. Sibony, Smooth
pseudoconvex domains in $\CC^2$ for which the corona theorem
and $L^p$ estimates for $\overline\partial$ fail, {\it Complex
analysis and geometry}, 209--222, Univ. Ser. Math., Plenum,
New York, 1993.

\item[{\bf [GAM]}]  T. W. Gamelin, {\it Uniform Algebras}, Prentice-Hall,
Englewood Cliffs, NJ, 1969.

\item[{\bf [GAR]}]  J. B. Garnett, {\it Bounded Analytic Functions}, 
Academic Press, New York, 1981.

\item[{\bf [GAJ]}] J. B. Garnett and P. W. Jones, The corona
theorem for Denjoy domains, {\it Acta Math.} 155(1985),
27--40.

\item[{\bf [GEL1]}]  I. M. Gelfand, To the theory of normed rings.  II.  On absolutely
convergent trigonometrical series and integrals, {\it C. R. (Doklady) Acad.\
Sci.\ URSS (N.S.)} 25(1939), 570--572.

\item[{\bf [GEL2]}]  I. M. Gelfand, To the theory of normed rings.  III.  On the ring
of almost periodic functions, {\it C. R. (Doklady) Acad.\
Sci.\ URSS (N.S.)} 25(1939), 573--574.

\item[{\bf [HAS]}] M. Hakim and N. Sibony, Fonctions
holomorphes born\'{e}es et limites tangentielles, {\it Duke
Math.\ Journal} 50(1983), 133--141.

\item[{\bf [KRA1]}]  S. G. Krantz, {\it Function Theory of Several Complex Variables},
$2^{\rm nd}$ ed., American Mathematical Society, Providence, RI, 2001.

\item[{\bf [KRA2]}] S. G. Krantz, {\it Cornerstones of
Geometric Function Theory: Explorations in Complex Analysis},
Birkh\"{a}user Publishing, Boston, 2006.

\item[{\bf [KRL1]}]  S. G. Krantz and S. Y. Li, Some remarks on
the corona problem on strongly pseudoconvex domains in
$\CC^n$, {\it Illinois J. Math.} 39(1995), 323--349.

\item[{\bf [KRL2]}] S. G. Krantz and S. Y. Li, Explicit
solutions for the corona problem with Lipschitz data in the
polydisc, {\it Pacific J. Math.} 174(1996), 443--458.

\item[{\bf [KRL3]}] S. G. Krantz and S. Y. Li, Factorization of
functions in subspaces of $L^1$ and applications to the corona
problem, {\it Indiana Univ.\ Math.\ J.} 45(1996), 83--102.

\item[{\bf [SIB]}] N. Sibony, Probl\`{e}me de la couronne pour
des domaines pseudoconvexes \`{a} bord lisse, {\it Ann.\ of
Math.} 126(1987), 675--682.

\end{enumerate}
\vspace*{.17in}

\begin{quote}
Department of Mathematics \\
Washington University in St.\ Louis  \\
St.\ Louis, Missouri 63130 \ \ U.S.A.  \\
{\tt sk@math.wustl.edu}
\end{quote}

\end{document}